\documentclass[10pt,draft]{article}
\usepackage{amssymb,amsmath}
\usepackage[english]{babel}
\usepackage[ansinew]{inputenc}
\usepackage{color}
\newcommand{\changeOK}[1]{{\color{blue}{#1}}}
\newcommand{\changesran}[1]{{\color{blue}{#1}}}
\def\normal{\renewcommand{\arraystretch}{1}}
\def\obec#1{{\normal \begin{tabular}[t]{c} #1 \end{tabular}}}

\newtheorem{theo}{Theorem}[section]
\newtheorem{lema}[theo]{Lemma}%[section]
\newtheorem{prop}{Proposition}[section]

\newcommand{\racion}[2]{\mbox{\small$\frac{{#1}}{{#2}}$}}
\newcommand{\half}{\frac{1}{2}}

\newcommand\be{\begin{enumerate}}
\newcommand\ee{\end{enumerate}}
\newcommand\bi{\begin{itemize}}
\newcommand\ei{\end{itemize}}

\newcommand{\pe}[2]{\langle#1,#2\rangle}

\newcommand\bd{\begin{definicion}{\bf }}
\newcommand\ed{\end{definicion}}
\newcommand\bl{\begin{lema}{\bf }}
\newcommand\el{\end{lema}}
\newcommand\bp{\begin{prop}{\bf }}
\newcommand\ep{\end{prop}}
\newcommand\bt{\begin{theo}{\bf }}
\newcommand\et{\end{theo}}
\newcommand\bdm{\begin{proof}}
\newcommand\edm{\end{proof}}
\newcommand\bn{\begin{nota}{\bf }}
\newcommand\en{\end{nota}}
\newcommand\bc{\begin{corolary}{\bf }}
\newcommand\ec{\end{corolary}}

\newtheorem{nota}{Remark}[section]

\newtheorem{definicion}{Definition}[section]
\newtheorem{corolary}{Corolary}[section]
\newenvironment{proof}{\noindent{\bf Proof:\ }}{\hfill
\fbox \par\vspace{2ex}}

\newfont{\got}{eufm10 scaled \magstep1}

\title{Limit relations between $q$-Krall type orthogonal
polynomials}
\date{\today}

\author{R. \'Alvarez-Nodarse${}^\dag{}^\ddag$ and R. S.
Costas-Santos${}^*$\\[5mm]
\small ${}^\dag$ Departamento de An\'alisis Matem\'atico.\\
\small Universidad de Sevilla. Apdo. 1160, E-41080 Sevilla, Spain\\
\small${}^\ddag$ Instituto Carlos I de F\'{\i}sica Te\'orica y
Computacional, \\
\small Universidad de Granada, E-18071 Granada, Spain\\
\small${}^*$ Departamento de Matem\'aticas,
E.P.S., Universidad Carlos III de Madrid.\\
\small Ave. Universidad 30, E-28911, Legan\'es, Madrid, Spain }
\begin{document}
\maketitle
%%%%%%%%%%%%%%%%%%%%%%%%%%%%%%%%%%%%%%%%%%%%%%%%%%%%%%%%%%%%%%%%%%%
\begin{abstract} \noindent
\small In this paper, we consider a natural extension of several
results related to Krall-type polynomials  introducing  a
modification of a $q$-classical linear functional via the addition
of one or two mass points.
The limit relations between the $q$-Krall type modification of big
$q$-Jacobi, little $q$-Jacobi, big $q$-Laguerre, and other families
of the $q$-Hahn tableau are established.
\end{abstract}
%%%%%%%%%%%%%%%%%%%%%%%%%%%%%%%%%%%%%%%%%%%%%%%%%%%%%%%%%%%%%%%%%%%
\section{Introduction}
In the last years, perturbations of a linear functional ${\cal C}$
via the addition of Dirac delta functions ---the so-called
Krall-type orthogonal polynomials--- have been intensively studied
(for recent reviews see e.g. \cite{ran-mar3,ran-mar-pet} and
references therein), i.e. ${\cal U}={\cal C}+A\delta(x-x_0)$,
where $A \ge0$, $x_0 \in \mathbb{R}$ and $\delta(x-y)$
means the Dirac linear functional defined by $\pe {\delta(x-y)}
{p(x)}=p(y)$, $\forall p\in \mathbb{C}[x]$, the linear space of
polynomials with complex coefficients.
Of particular interest are the cases when the starting functional
is a classical linear functional (Jacobi \cite{ran-mar3,koo84},
Laguerre \cite{ran-mar3,koko}, Hermite \cite{ran-mar3}, and Bessel
\cite{mar-mar}) and a discrete one (Hahn, Meixner, Kravchuk, and
Charlier) \cite{ran-gar-mar,ran-mar1,ran-mar2,hbh,bavko,god97}.
A more general case ${\cal U}={\cal C}+\sum_{i=1}^M A_i\delta
(x-a_i)-\sum_{j=1}^NB_j\delta'(x-b_j)$ was studied in a recent
paper \cite{ran-arv-mar} where a special emphasis is given when
${\cal C}$ is a semiclassical linear functional.

In a recent paper \cite{ran-pet} the case when  $\cal C$ is a
discrete semiclassical or $q$-semiclassical linear functional was
considered in details.
Here we will focus our attention on the case when $\cal C$ is a
$q$-classical linear functional and we will construct the
Krall-type polynomials associated with the $q$-classical families
of the so-called $q$-Hahn tableau \cite{ran01,koo94}.
This case is not so well known and only few papers deals with
examples of such polynomials: the Stieltjes-Wigert polynomials
\cite{chi1}, a particular case of the $q$-little Jacobi
polynomials \cite{vin01}, and the Al-Salam \& Carlitz I and
discrete $q$-Hermite I \cite{ran-pet}.

The aim of the present contribution\footnote{
\color{red}This paper appear in R. \'Alvarez-Nodarse and R.S. Costas-Santos,
Limit relations between q-Krall type orthogonal polynomials, 
{\it J. Math. Anal. Appl.} {\bf 322}(1) (2006) 158-176, but we have found a misprint in
formulas \eqref{ker-rep1} and  \eqref{ker-rep2}. In this version we correct the misprints
as well as all all formulas that changes as a consequence of such misprint. The 
changes are in blue.} is to continue the work
started in \cite{ran-pet} and study several families of $q$-Krall
type orthogonal polynomials.
In particular, we will obtain the limits of the $q$-Krall type
polynomials in the $q$-Hahn tableau.
In such a way we will continue the study started in \cite{ran-mar3}
concerning the limit relations among the Krall-type families.

The structure of the paper is the following. In section 2 some
preliminaries and the basic parameters of the families that we will
consider later on are given. In particular, we include the
explicit values for the kernels of the corresponding $q$-classical
polynomials in terms of the polynomials and their $q$-derivatives.
In section 3 the $q$-Krall-type orthogonal polynomials are defined
and some algebraic properties are deduced for these new families.
Finally, in section 4, the limits of the modified polynomials of
the examples considered in section 2 are established.
%%%%%%%%%%%%%%%%%%%%%%%%%%%%%%%%%%%%%%%%%%%%%%%%%%%%%%%%%%%%%%%%%%%
\section{Preliminary results}
In this section, we state some formulas for $q$-classical
orthogonal polynomials $P_n(x(s))_q\equiv P_n(s)_q$ of the
$q$-Hahn tableau,
orthogonal with respect to a $q$-classical linear functional
${\cal C}_q$ \cite{ran02}, i.e.,
\begin{equation}\label{tab-func-C}
\pe {{\cal C}_q} {P_n P_m}=d_n^2 \delta_{n,m}, \quad d_n^2\ne 0,
\quad n, m= 0, 1, 2, \dots
\end{equation}
These functionals usually have the form (see the section
\ref{sect-examp} for more details)
\[
\pe {{\cal C}_q} P = \left\{ \begin{array}{ll} \displaystyle
\sum_{s=0}^\infty P(x(s)) \rho(s) \nabla x_1(s), & \ \mbox{\small
little $q$-Jacobi, $q$-Meixner, Wall, $q$-Charlier,}
\\[5mm] \displaystyle \int_{s_0}^{s_1} P(x) \rho(x)\,  d_qx, &
\ \mbox{\small big $q$-Jacobi, big $q$-Laguerre, Al-Salam-Carlitz
I,} \end{array} \right.
\]
etc., where $\int_{s_0}^{s_1}f(t)d_qt$ is the Jackson $q$-integral
(see \cite{gasper1,ks}), $\rho$ is a weight function satisfying the
following difference equation of Pearson-type
$$
\Delta [\sigma(s) \rho(s)] =\tau(s)\rho(s) \nabla x_1(s) \ \iff \
\frac{\rho(s+1)}{\rho(s)}= \frac{\sigma(s)+\tau(s)\nabla x_1(s)}
{\sigma(s+1)},
$$
the lattice is $x(s)=c q^{\pm s}+c'$, $x_k(s)=x(s+\racion k 2)$,
and $\nabla$ and $\Delta$ are the backward and forward difference
operators defined respectively  as $\nabla f(s)=f(s)-f(s-1)$,
$\Delta f(s)=f(s+1)-f(s)$.
Now, consider the sequence of $q$-classical orthogonal polynomials
with respect to the linear functional ${\cal C}_q$ ($q$-COP).
They satisfy the second order linear difference equation (SODE) of
hypergeometric type \cite{nsu} $\sigma(s) \frac{\Delta }{\nabla
x_1(s)} \frac{\nabla y(s)}{\nabla x(s)}+\tau(s) \frac{\Delta y(s)}
{\Delta x(s)} +\lambda_n y(s)=0$, where $\sigma(s)$ and $\tau(s)$
are polynomials of degree at most 2 and exactly 1, respectively,
and $\lambda_n$ is a constant. Moreover, these families of
$q$-polynomials satisfy several algebraic relations such as a
three-term recurrence relation (TTRR)
\begin{equation} \label{q-TTRR}
x(s)P_n(s)_q=\alpha_n P_{n+1}(s)_q+\beta_n
P_n(s)_q+\gamma_n P_{n-1}(s)_q, \  \ n=0,1,2\dots,
\end{equation}
with the initial conditions $P_0(s)_q=1$, $P_{-1}(s)_q=0$, the
structure relations ($n=1,2,3,\dots$)
\begin{equation} \label{rel-str} \begin{split} \displaystyle
\sigma(s)\frac{\nabla P_n(s)_q}{\nabla x(s)} & = \widetilde\alpha_n
P_{n+1}(s)_q+\widetilde\beta_n P_n(s)_q+\widetilde \gamma_n
P_{n-1}(s)_q, \\[3mm] \displaystyle \phi(s)
\frac{\Delta P_n(s)_q}{\Delta x(s)}  & =  \widehat\alpha_n
P_{n+1}(s)_q + \widehat\beta_n P_n(s)_q+\widehat \gamma_n
P_{n-1}(s)_q, \end{split}
\end{equation}
where $\phi(s)=\sigma(s)+\tau(s) \nabla x_1(s)$, as well as the
Christoffel-Darboux formula for the $n$-th kernel associated with
the family
\begin{equation} \label{chr-dar}
K_{n}(s_1,s_2):=\sum_{m=0}^n \frac{P_m(s_1)_q P_m(s_2)_q}{d_m^2}=
\frac{\alpha_n}{d_{n}^2} \frac{P_{n+1}(s_1)_q P_{n}(s_2)_q -
P_{n+1}(s_2)_q P_{n}(s_1)_q}{x(s_1)-x(s_2)}.
\end{equation}
In the sequel we will use the notation
$\mathbb{K}_{n}(s_0):=K_n(s_0,s_0)$.
From \eqref{chr-dar} and  \eqref{rel-str} follows that
\be
\item If $\sigma(s_0)=0$, then
{
\begin{equation} \label{ker-rep1} \displaystyle {\color{blue} K_{n}(s,s_0)}=
\frac{P_n(s_0)_q}{d_{n}^2\left(\racion{\widetilde \gamma_n}
{\gamma_n}-\racion{\widetilde \alpha_n}{\alpha_n}\right)}\left[
\frac {\widetilde \gamma_n}{\gamma_n} P_n(s)_q-\frac{\sigma(s)}
{x(s)-x(s_0)} \frac{\nabla P_n(s)_q}{\nabla x(s)}\right].
\end{equation}
}
\item If $\phi(s_0)=0$, then
{
\begin{equation} \label{ker-rep2} {\color{blue} K_{n}(s,s_0)}= \displaystyle
\frac{P_{n}(s_0)_q}{d_{n}^2\left(\racion {\widehat
\gamma_n}{\gamma_n}- \racion{\widehat
\alpha_n}{\alpha_n}\right)}\left[\frac{\widehat \gamma_n}{\gamma_n}
P_n(s)_q- \frac{\phi(s)}{x(s)-x(s_0)} \frac{\Delta P_n(s)_q}{\Delta
x(s)} \right].
\end{equation}}
\ee
\bn A straightforward calculation shows that $\displaystyle
\frac {\widetilde \alpha_n} {\alpha_n}$  and $\displaystyle
\frac{\widetilde \gamma_n}{\gamma_n}$ are independent of the
normalization of $P_n(s)_q$, i.e. if $\widehat P_n(s)=C_n P_n(s)_q$
then those ratios do not change.
Moreover (\cite[eq. (6.15)]{ran1}) ${\widetilde \gamma_n}/
{\gamma_n}-{\widetilde \alpha_n}/{\alpha_n}$ $=$${\widehat
\gamma_n}/{\gamma_n}-{\widehat \alpha_n}/{\alpha_n}$.
\en
%%%%%%%%%%%%%%%%%%%%%%%%%%%%%%%%%%%%%%%%%%%%%%%%%%%%%%%%%%%%%%%%%%%
\subsection{The $q$-classical polynomials} \label{sect-examp}
In this section, we will summarize the main properties of the
$q$-polynomials of the $q$-Hahn tableau needed in the next sections
(for more details see \cite{ks}). In all cases we have used
\eqref{ker-rep1} and \eqref{ker-rep2} for computing the kernels at
the corresponding points.
In the sequel we will consider probabilistic measures, i.e.
$d_0^2=1$.
This fact will be useful in order to obtain the right limits for
the corresponding $q$-Krall-type polynomials. Here and through out
the paper we will use the standard notation for the basic series.
For more details see \cite{gasper1}.
%%%%%%%%%%%%%%%%%%%%%%%%%%%%%%%%%%%%%%%%%%%%%%%%%%%%%%%%%%%%%%%%%%%
\begin{table}[t]
\begin{footnotesize}
\caption{Parameters of big $q$-Jacobi and Stieltjes-Wigert
polynomials\label{t1}}
\begin{center}
\begin{tabular}{l@{\hspace{27mm}}l@{\hspace{27mm}}l}
\hline
& Big $q$-Jacobi & Stieltjes-Wigert \\
\hline $\sigma(x)$ & $q^{-1}(x-aq)(x-cq)$ & $q^{-1}x$\\
$\phi(s)$ &  $aq(x-1)(bx-c)$ & $x^2$\\ $\lambda_n$ & $-q^{\half}
\frac{(1-abq^{n+1})(1-q^n)} {(1-q)^2q^n}$ & $q^\half\frac{1-q^n}
{(1-q)^2}$ \\ $\rho(x)$ & $\frac{1}{aq(1-q)}\frac{(aq,bq,cq,
abc^{-1}q,a^{-1}x, c^{-1}x;q)_\infty}{(q,abq^2,a^{-1}c,ac^{-1}q,x,
bc^{-1}x;q)_{\infty}}$ & $-\frac{1}{\ln q}\frac 1{(q,-x,-qx^{-1};
q)_\infty}$ \\[3mm] $d_n^2$ & $\frac{(1-abq)(q,bq,abc^{-1}q;q)_n}
{(1-abq^{2n+1})(aq,abq,cq;q)_n}\left(-acq^{\frac {n+3}2}\right)^n$
& $\frac 1 {(q;q)_nq^{n}}$\\[3mm] $P_n(x_0)$ & $P_n(aq;a,b,c;q)=
\frac{(abc^{-1}q;q)_n}{(cq;q)_n}\left(-cq^{\frac{n+1} 2}\right)^n$
& $S_n(0;q)=\frac 1 {(q;q)_n}$\\[3mm] $P_n(x_1)$ & $P_n(cq;a,b,c;
q)=\frac{(bq;q)_n}{(aq;q)_n} \left(-aq^{\frac {n+1} 2}\right)^n$ &
\mbox{---} \\ ${\widehat \alpha_n}/{\alpha_n}$ & $abq\frac {1-q^n}
{1-q}$ & $\frac{1-q^n}{1-q}$ \\ ${\widehat \gamma_n}/{\gamma_n}$ &
$-\frac {1-abq^{n+1}}{(1-q)q^n}$ &  $\frac{1}{1-q}$\\ ${\widetilde
\alpha_n}/{\alpha_n}$ & $\frac {1-q^n}{(1-q)q^n}$ & $0$\\
${\widetilde \gamma_n}/{\gamma_n}$ & $-\frac {1-abq^{n+1}}{1-q}$ &
$\frac {q^n}{1-q}$\\ \hline
\end{tabular}
\end{center}
\end{footnotesize}
\end{table}

The big $q$-Jacobi polynomials $P_n(x;a,b,c;q)$, introduced by Hahn
in 1949, are the most general family of $q$-polynomials on the
$q$-linear lattice $x:=x(s)=q^s$.
They constitute a $q$-COP sequence with respect to the linear
functional ${\cal C}^{BqJ}$
\begin{equation}\label{wf-bqJ}
\pe {{\cal C}^{BqJ}} P:=\int_{cq}^{aq} P(x) \rho(x)d_q x,
\end{equation}
where the weight function $\rho(x)$, supported on $[cq,aq]$,
$0<a,b<q^{-1}, c<0$, is given in table \ref{t1}.

For these polynomials we also need the following expressions for
the kernels
{\begin{small}
\begin{equation} \label{BqJ-3}
\begin{split}
{\color{blue} K_{n}^{BqJ}(x,aq)}=& \displaystyle \frac{(aq,abq;q)_n}{(q,bq;q)_n}
\left[\frac{(1-abq^{n+1})P_n(x;a,b,c;q)-(x-cq)(1-q^{-1}) {\cal
D}_{q} P_n(q^{-1}x;a,b,c;q)}{(1-abq)a^n}\right], \\[3mm]
{\color{blue} K_{n}^{BqJ}(x,1)} =&  \displaystyle \frac{(aq,abq,cq;q)_n}
{(q,bq,abc^{-1}q;q)_n}\left[\frac{(q^{-n}-abq)P_n(x;a,b,c;q)\!-\!
aq(bx\!-\!c)(q\!-\!1){\cal D}_q P_n(x;a,b,c;q)}{(1-abq)
(-acq^{n/2+1/2})^n}\right],\\[3mm] 
{\color{blue} K_{n}^{BqJ}(x,cq)} = &
\displaystyle \frac{(abq,cq;q)_n}{(q,abc^{-1}q;q)_n}\left[\frac
{(1-abq^{n+1})P_n(x;a,b,c;q)\!-\!(x\!-\!aq)(1\!-\!q^{-1}){\cal
D}_{q} P_n(q^{-1}x;a,b,c;q)}{(1-abq)c^n}\right], \end{split}
\end{equation}
\end{small}}%
and
{
\[
\begin{array}{rl}
{\color{blue} \mathbb{K}_{n}^{BqJ}(aq)}= & \displaystyle {\color{blue} \sum_{k=0}^{n}}
\frac{(1-abq^{2k+1})(aq,abq,abc^{-1}q;q)_k}{(1-abq)(q,bq,cq;q)_k}
\left(-a^{-1}cq^{\racion {k-1}2}\right)^{k},\\[4mm]
{\color{blue} \mathbb{K}_{n}^{BqJ}(cq)}= & \displaystyle {\color{blue} \sum_{k=0}^{n}}
\frac{(1-abq^{2k+1})(bq,abq,cq;q)_k}{(1-abq)(q,aq,abc^{-1}q;q)_k}
\left(-ac^{-1}q^{\racion {k-1}2}\right)^{k},\\[4mm]
{\color{blue} \mathbb{K}_{n}^{BqJ}(1)}=& \displaystyle  {\color{blue} \sum_{k=0}^{n}}
\frac{(1-abq^{2k+1})(aq,abq,cq;q)_k}{(1-abq)(q,bq,abc^{-1}q;q)_k}
\left(-acq^{\racion {k+3}2}\right)^{-k},\\ 
{\color{blue} K_{n}^{BqJ}(aq,cq)}=
& \displaystyle {\color{blue} \sum_{k=0}^{n} } \frac{(1-abq^{2k+1})(abq;q)_k
\left(-q^{\frac{k-1} 2}\right)^k}{(1-abq)(q;q)_k}={\color{blue} \frac{(abq^2;
q)_{n}}{(q;q)_{n}} \left(-q^{\racion {n+1}2}\right)^{n}},
\\[4mm] 
{\color{blue} K_{n}^{BqJ}(aq,1)}=& \displaystyle  {\color{blue} \sum_{k=0}^{n}}
\frac{(1-abq^{2k+1})(aq,abq;q)_k}{(1-abq)(q,bq;q)_k} (aq)^{-k}={\color{blue} 
\frac{(aq^2,abq^2;q)_{n}}{(q,bq;q)_{n}}(aq)^{-n}},
\end{array}
\]}
where ${\cal D}_q$ the $q$-Jackson derivative (see e.g. \cite{ks}),
${\cal D}_q P(z)=[P(z)-P(qz)]/[(1-q)z]$.

The big $q$-Laguerre polynomials $P_n(x;a,c;q)$, are a particular
case of the big $q$-Jacobi polynomials: $P_n(x;a,c;q)=P_n(x;a,0,c;
q)$, therefore all their properties can be obtained from the
corresponding ones of the big $q$-Jacobi by putting $b=0$.
A special case of the big $q$-Laguerre polynomials are the
affine $q$-Kravchuk polynomials \cite[page 101]{ks}.

\begin{table}[t]
\begin{footnotesize}
\caption{Parameters of little $q$-Jacobi, $q$-Laguerre and
Al-Salam-Carlitz I polynomials\label{t2}}
\begin{center}\begin{tabular}
{l@{\hspace{11mm}}l@{\hspace{11mm}}l@{\hspace{11mm}}l}
\hline
& little $q$-Jacobi & $q$-Laguerre & Al-Salam-Carlitz I\\
\hline
$\sigma(x)$ & $q^{-1}x(x-1)$ & $q^{-1}x$ & $(x-1)(x-a)$\\
$\phi(s)$ &  $ax(bqx-1)$ & $ax(x+1)$ & $a$ \\ $\lambda_n$ &
$-q^{\frac 3 2}\frac{(1-abq^{n+1})(1-q^n)}{(1-q)^2q^n}$ &
$aq^\half\frac{1-q^n}{(1-q)^2}$ & $-q^{\frac 3 2}\frac{1-q^n}
{(1-q)^2q^n}$\\ $\rho(x)$ & $\frac{(aq;q)_\infty}{(abq^2;q)_\infty}
\frac{(bq;q)_s}{(q;q)_s}\, a^s$ & $\frac{ a^s\, (aq,-c,-c^{-1}q;
q)_\infty}{(q,-caq,-c^{-1}a^{-1},-cq^s;q)_\infty}$ & $\frac{1}{1-q}
\frac{(qx,a^{-1}qx;q)_\infty}{(q,a,a^{-1}q;q)_\infty}$\\[3mm]
$d_n^2$ &  $\frac{(1-abq)}{(1-abq^{2n+1})}\frac{(q,bq;q)_n}
{(aq,abq;q)_n}(aq)^n$ & $\frac{(aq;q)_n}{(q;q)_n q^n}$ & $(q;q)_n
\left(-a q^{\frac {n-1} 2}\right)^n$  \\[3mm] $P_n(x_0)$ &
$p_n(0;a,b|q)=1$ & $L_n^{(\alpha)}(0;q)=\frac{(aq;q)_n}{(q;q)_n}$ &
$U_n^{(a)}(1;q)=\left(-aq^{\frac {n-1}2}\right)^n$\\[3mm]
${\widetilde \alpha_n}/{\alpha_n}$ & $\frac{1-q^n}{(1-q)q^n}$ &
$0$ & $q\frac{(1-q^n)}{(1-q)q^n}$\\ ${\widetilde\gamma_n}/
{\gamma_n}$ & $-\frac{1-abq^{n+1}}{1-q}$ & $\frac {aq^n}{1-q}$ &
$-\frac {q}{1-q}$\\  \hline
\end{tabular}\end{center}
\end{footnotesize}
\end{table}

The little $q$-Jacobi polynomials $p_n(x;a,b|q)$ constitute a
$q$-OPS with respect to a linear functional
$$
\pe {{\cal C}^{lqJ}} P :=\sum_{s=0}^\infty P(s) \rho(s) q^s,
$$
where $\rho(s)$ is given in table \ref{t2} and it is supported on
$[0,1]$, $0<a<q^{-1}, \ b<q^{-1}$.
Moreover
{
\begin{equation} \label{lqJ}
\begin{array}{c}
\displaystyle {\color{blue} K_{n}^{lqJ}(x,0)}=
\frac{(aq,abq;q)_n}{(q,bq;q)_n a^n}\left[\frac{(1-abq^{n+1})
P_n(x;a,b|q)-(x-1)(1-q^{-1}) {\cal D}_q P_n(q^{-1}x;a,b|q)}{1-abq}
\right], \\[4mm] 
\displaystyle {\color{blue} \mathbb{K}_{n}^{lqJ}(0)}=
\sum_{k=0}^{n}\frac{(1-abq^{2k+1})(aq,abq;q)_k}
{(1-abq)(q,bq;q)_k}\, (aq)^{-k}=\frac{(aq^2,abq^2;q)_{n}}
{(q,bq;q)_{n}}(aq)^{-n}.
\end{array}
\end{equation}}

The $q$-Meixner polynomials $M_n(q^{-s};b,c;q)$ are a $q$-COP
sequence with respect to a linear functional %${\cal C}^{qM}$
$$
\pe {{\cal C}^{qM}} P:=\sum_{s=0}^\infty  P(s)\rho(s)q^{-s},
$$
where the weight function $\rho(s)$ is supported on $[1,+\infty)$,
$0<b<q^{-1}, \ 0<c$ (see table \ref{t3}).
Furthermore,
{
\begin{equation} \label{M}
\begin{array}{c}
\displaystyle {\color{blue} K_{n}^{qM}(x,1)}=\frac{(bq;q)_n}{(q,-c^{-1}q;q)_n}
\left[M_n(x;b,c;q)-(x+bc)(1-q){\cal D}_{q^{-1}} M_n(x;b,c;q)
\right],\\[4mm] \displaystyle 
{\color{blue} \mathbb{K}_{n}^{qM}(1)}=
{\color{blue} \sum_{k=0}^{n}}\frac{(bq;q)_k}{(q,-c^{-1}q;q)_k}\, q^k.
\end{array}
\end{equation}}
A special case of the $q$-Meixner polynomials are the quantum
$q$-Kravchuk \cite[page 98]{ks}.

\begin{table}[t]
\begin{footnotesize}
\caption{Parameters of Wall, $q$-Meixner, and $q$-Charlier
polynomials\label{t3}}
\begin{center}\begin{tabular}
{l@{\hspace{14mm}}l@{\hspace{14mm}}l@{\hspace{14mm}}l}
\hline & Wall & $q$-Meixner & $q$-Charlier \\ \hline $\sigma(x)$ &
$q^{-1}x(x-1)$ & $q^{-1}c(x-bq)$ & $q^{-1}ax$ \\ $\phi(s)$ & $-ax$&
$(x-1)(x+bc)$ & $x(x-1)$\\ $\lambda_n$ & $-q^{\frac 32}\frac{1-q^n}
{(1-q)^2q^n}$ &$-q^{\half}\frac{1-q^n}{(1-q)^2q^n}$ & $q^\half
\frac{1-q^n}{(1-q)^2}$\\ $\rho(x)$ & $\frac{(aq;q)_\infty}{(q;q)_s}
a^s$ & $\frac{(-bcq;q)_\infty}{(-c;q)_\infty}\frac{(bq;q)_s}
{(q,-bcq;q)_s}\, \left(c q^{\frac {s+1} 2}\right)^s$ & $
\frac{q^{-\half}} {(q;q)_s (-a;q)_\infty} \, \left(aq^{\racion{s+1}
2}\right)^s$ \\[3mm] $d_n^2$ &  $\frac{(q;q)_n}{(aq;q)_n} (aq)^n$ &
$\frac{(q,-c^{-1}q;q)_n}{(bq;q)_n}q^{-n}$ &
$(-a^{-1}q,q;q)_n\, q^{-n}$ \\[3mm] $P_n(x_0)$ & $p_n(0;a|q)=1$ &
$M_n(1;b,c|q)=1$ & $C_n(1;a;q)=1$ \\[3mm] ${\widehat \alpha_n}/
{\alpha_n}$ & $0$ & $\frac{1-q^n}{1-q}$ & $\frac{1-q^n}{1-q}$ \\
${\widehat \gamma_n}/ {\gamma_n}$ & $-\frac {1}{(1-q)q^n}$ &
$\frac {1}{1-q}$ & $\frac {1}{1-q}$ \\ \hline
\end{tabular}\end{center}
\end{footnotesize}
\end{table}

The Al-Salam-Carlitz I polynomials $U_n^{(a)}(x;q)$ are orthogonal
with respect to the linear functional %${\cal C}^{ACI}$
$$
\pe {{\cal C}^{ACI}} P:=\int_a^1  P(x) \rho(x) d_q x,
$$
where $\rho(x)$  is supported on $[a,1]$, $a<0$, $x:=x(s)=q^s$.
Their main data are in table \ref{t3}. For these polynomials we
have
{
\begin{equation} \label{ACI}
\begin{array}{c}
\displaystyle{\color{blue}  K_{n}^{ACI}(x,1)}=\frac{q^n}{(q;q)_n}\left[U_n^{(a)}
(x;q)-(x-a)(1-q^{-1}){\cal D}_q U_n^{(a)}(q^{-1}x;q)\right],\\[4mm]
\displaystyle {\color{blue} \mathbb{K}_{n}^{ACI}(1)=\sum_{k=0}^{n}} \frac{1}
{(q;q)_k}\left(-aq^{\racion {k-1} 2}\right)^k.
\end{array}
\end{equation}}

The little $q$-Laguerre / Wall polynomials $p_n(x;a|q)$ are
orthogonal with respect to the linear functional %${\cal C}^{lqL}$
$$
\pe {{\cal C}^{lqL}} P:= \sum_{s=0}^\infty P(s) \rho(s)q^s,\quad
x:=x(s)=q^s,\, \mbox{supp}(\rho)=[0,1].
$$
Since they are a particular case of little $q$-Jacobi ($b=0$) all
their properties can be obtained from the former ones putting
$b=0$ (see table \ref{t3}).  In particular,
\changeOK{\begin{equation} \label{lqL}
\begin{array}{c}
\displaystyle K_{n}^{lqL}(x,0)=
\frac{ (aq;q)_n}{(q;q)_na^{n}q^n}\left[p_n(x;a|q)
-a(1-q)q^n{\cal D}_q p_n(x;a|q)\right],\quad
\displaystyle \mathbb{K}_{n}^{lqL}(0)= \frac{(aq^2;q)_{n}}
{(q;q)_{n}(aq)^{n}}.
\end{array}
\end{equation}}

The $q$-Laguerre polynomials $L^{(\alpha)}_n(x;q)$ are orthogonal
with respect to the linear functional %$\mathcal{C}^{qL}$
\[
\pe{\mathcal{C}^{qL}}{P}:= \sum_{s=-\infty}^{+\infty} P(cq^k)
\rho(s)q^s,
\]
where the weight function $\rho(s)$ (see table \ref{t2}) is
supported on $[0,+\infty)$, $a=q^\alpha$, $x:=x(s)=cq^s$.
In this case
{\begin{equation} \label{L}
\begin{array}{c} 
{\color{blue} \displaystyle K_{n}^{qL}(x,0)=q^nL_n^{(\alpha)}
(x;q)-\frac{q^{-1}-1} {a}{\cal D}_q L_n^{(\alpha)}(q^{-1}x;q),\quad
\displaystyle \mathbb{K}_{n}^{qL}(0)=\sum_{k=0}^{n} \frac{(aq;q)_k}
{(q;q)_k}\,q^k= \frac{(aq^2;q)_{n}}{(q;q)_{n}}. }\end{array}
\end{equation}}
The $q$-Charlier polynomials $C_n(q^{-s};a;q)$ constitute a
$q$-COP sequence with respect to the linear functional
%${\cal C}^{qC}$
$$
\pe {{\cal C}^{qC}} P:=\sum_{s=0}^\infty  P(s) \rho(s)q^{-s},
$$
where $\rho(s)$ is supported on $[1,+\infty)$, $a>0$ (see table
\ref{t3}). Moreover,
{
\begin{equation} \label{C}
\begin{array}{c} \displaystyle 
{\color{blue} K_{n}^{qC}(x,1)=\frac{C_n(x;a;q)-
x(1-q){\cal D}_{q^{-1}} C_n(x;a;q)}{(-a^{-1}q,q;q)_n},\quad
\displaystyle\mathbb{K}_{n}^{qC}(1)=\sum_{k=0}^n \frac{q^k}
{(-a^{-1}q,q;q)_k}. }
\end{array}
\end{equation}}
The Stieltjes-Wigert polynomials $S_n(x;q)$ correspond to an
indeterminate moment problem, so there are infinitely many
representations for the linear functional ${\cal C}^{SW}$ with
respect to which they are orthogonal (see e.g. \cite{ks}).
Here we will chose the following one
$$
\pe {{\cal C}^{SW}}P :=\int_0^\infty P(x)\rho(x)dx,
$$
where $\rho(s)$ is a weight function supported on $[0,+\infty)$
(see table \ref{t1}). In this case
{
\begin{equation} \label{SW}
\begin{array}{c}
\displaystyle {\color{blue} K_{n}^{SW}(x,0)=q^nS_n(x;q)-(q^{-1}-1){\cal D}_q
S_n(q^{-1}x;q), \quad \displaystyle \mathbb{K}_{n}^{SW}(0)=\frac
1{(q;q)_{n}}. }
\end{array}
\end{equation}}
%%%%%%%%%%%%%%%%%%%%%%%%%%%%%%%%%%%%%%%%%%%%%%%%%%%%%%%%%%%%%%%%%%%
\section{The $q$-Krall-type orthogonal polynomials}
In this section, we will introduce the $q$-Krall-type
orthogonal polynomials. In a very recent paper \cite{ran-pet}
the authors introduce the ``discrete'' Krall polynomials as a
perturbation of a classical or semiclassical discrete linear
functional and they develop a general theory in order to find some
algebraic properties such as TTRR, SODE, etc.
In this paper we focus our attention on the special case when the
starting functional $\cal C$ is a $q$-classical functional
\cite{ran02}. Thus we consider the linear functional ${\cal U}$
defined as
\begin{equation} \label{func-U}
\langle {\cal U},P\rangle = \langle {\cal C},P\rangle+A P(x_0)+
BP(x_1), \quad  A,\, B \ge 0,
\end{equation}
where $\cal C$ is the linear functional \eqref{tab-func-C} and
$x_0,\, x_1\in \mathbb{R}$.
In \cite{ran-pet} a general theory for solving this problem (when
$N$ mass points are added) has been presented, nevertheless only
two examples were considered in details.
Here we will complete this work introducing new examples and we
will establish the limit relation among them, in the same way as
in \cite{ran-mar3}.

The explicit expression of the  polynomials $\widetilde P_n^{A,B}
(s)_q$ orthogonal with respect to the linear functional ${\cal U}$
\eqref{func-U} is given by \cite{ran-pet} (it is assumed that
the polynomials $\widetilde P_n$ and $P_n$
have the same leading coefficient)
\begin{equation}\label{FR-ker}
\widetilde P_n(x)=P_n(x)-\sum_{i=1}^M A_i\widetilde P_n(a_i)
K_{n-1}(x,a_i),
\end{equation}
where $(\widetilde P_n(a_k))_{k=1}^M$ are the solution of the
system
$$
\widetilde P_n(a_k)= P_n(a_k)-\sum_{i=1}^M A_i\widetilde P_n(a_i)
K_{n-1}(a_k,a_i),\quad k=1,2,\dots, M.
$$
The formula \eqref{FR-ker} was firstly obtained by Uvarov
\cite{uva69} (see also \cite[\S 2.9]{ism06}).
Hence, the formula \cite[Eq. (2.5) page 57]{ran-pet} yields
\begin{equation} \label{gen-rep-pol-2-mas}
\widetilde P_n^{A,B}(s)_q=P_n(s)_q\!-\!\left[\!\begin{array}{c}
AK_{n-1}(x,x_0) \\ BK_{n-1}(x,x_1)\end{array} \!\right]^t
\!\left[\!\begin{array}{cc} 1+A\mathbb{K}_{n-1}(x_0) & B
K_{n-1}(x_0,x_1) \\ A K_{n-1}(x_1,x_0) &  1+B\mathbb{K}_{n-1}(x_1)
\end{array} \!\right]^{-1} \!\left[\!\begin{array}{c} P_{n}(x_0)
\\ P_{n}(x_1) \end{array} \!\right],
\end{equation}
where $C^t$ is the transpose of $C$.
Furthermore, the polynomials $\widetilde P_n^{A,B}(s)_q$  exist
for every $n=0,1,\dots$ if and only if the following condition
\begin{equation}\label{exi}
\det\left[\begin{array}{cc} 1+A\mathbb{K}_{n-1}(x_0) & B K_{n-1}
(x_0,x_1) \\ A K_{n-1}(x_1,x_0) &  1+B\mathbb{K}_{n-1}(x_1)
\end{array} \right]\ne 0, \quad \forall  n\in \mathbb{N},
\end{equation}
holds.
When the mass $B=0$  \eqref{gen-rep-pol-2-mas} transforms into
\begin{equation} \label{pol_one_mas}
\widetilde{P}_n^A(s)_q =
P_n(s)_q - A \widetilde{P}_n^A(x_0)_q K_{n-1} (x,x_0),\quad
\widetilde{P}_n^A(x_0)_q=
\frac{P_n(x_0)_q}{1+A \mathbb{K}_{n-1}(x_0)}.
\end{equation}
Notice that if $A\ge 0$, then \eqref{exi} becomes into
$1+A\mathbb{K}_{n-1}(x_0)\ge 1$, hence $\widetilde P_n^{A}(s)_q$
exists for every $n=0,1,\dots$.
\\
The next step is to construct the corresponding families of
$q$-Krall type orthogonal polynomials associated with each family
of $q$-orthogonal polynomials considered in section
\ref{sect-examp}. We will start with the  big $q$-Jacobi family
since the other families can be obtained from it via taking
appropriate limits.
Furthermore, we will choose the values of $x_0$ and $x_1$ in such
a way that the kernels \eqref{chr-dar} has the simplest form, i.e.,
\eqref{ker-rep1} and \eqref{ker-rep2}.
%%%%%%%%%%%%%%%%%%%%%%%%%%%%%%%%%%%%%%%%%%%%%%%%%%%%%%%%%%%%%%%%%%%
\subsection{The big $q$-Jacobi-Krall polynomials}
Let us consider the linear functional  ${\cal U}^{BqJ}$ defined by
$$
\pe {{\cal U}^{BqJ}}{P}=\pe{{\cal C}^{BqJ}}{P}+A P(x_0)+BP(x_1),
\quad A,\, B\ge 0,
$$
where $x_0,\, x_1\in\mathbb{R}$ and  ${\cal C}^{BqJ}$ is the
functional \eqref{wf-bqJ}.
The corresponding polynomials will be denoted by $\widetilde
P_n^{A,B}(x;a,b,c;q)$ and constitute a $q$-analog of the
Koornwinder polynomials \cite{koo84}.
The polynomial expression for this family follows from
\eqref{gen-rep-pol-2-mas}
\[
\begin{split}
\widetilde  P_n^{A,B}(x;a,b,c;q) = & P_n(x;a,b,c;q)-\left[
\begin{array}{cc}AK_{n-1}^{BqJ}(x;x_0) &  BK_{n-1}^{BqJ}(x;x_1)
\end{array} \right] \times\\ &  \left[\begin{array}{cc} 1+A
\mathbb{K}_{n-1}^{BqJ}(x_0) & BK_{n-1}^{BqJ}(x_0,x_1) \\
AK_{n-1}^{BqJ}(x_1,x_0) & 1+B\mathbb{K}_{n-1}^{BqJ}(x_1)
\end{array} \right]^{-1}\left[\begin{array}{c} P_n(x_0;a,b,c;q)
\\ P_n(x_1;a,b,c;q)\end{array}\right]. \end{split}
\]
Now, we are going to consider two specific cases:\bigskip

{\color{blue} 
\noindent {\bf 1.} The $q$-Koornwinder polynomials obtained when
we add two mass points at the endpoints of the interval of
orthogonality of the big $q$-Jacobi polynomials. i.e.,  $x_0=cq$
and $x_1=aq$.
For these values,
$$
\widetilde P^{A,B}_n(x;a,b,c;q):=P_n(x;a,b,c;q)-A\widetilde
P^{A,B}_n(cq) K^{BqJ}_{n-1}(x,cq)- B\widetilde
P^{A,B}_n(aq)K^{BqJ}_{n-1}(x,aq).
$$

Then, using \eqref{BqJ-3} and taking into account the identities
\cite[Eq. (3.5.6), (3.5.7)]{ks} for the big $q$-Jacobi polynomials,
\begin{eqnarray}\nonumber
P_n(x;a,b,c;q)-P_n(qx;a,b,c;q)=\frac{q(q^{-n}-1)(1-abq^{n+1})}
{(1-aq)(1-cq)}\,xP_{n-1}(qx;aq,bq,cq;q), \\ \label{rel-for-BqJ}
{\cal D}_q P_n(x;a,b,c;q)=\frac{q(q^{-n}-1)(1-abq^{n+1})}{(1-q)
(1-aq)(1-cq)}\, P_{n-1}(qx;aq,bq,cq;q),
\end{eqnarray}
we get
\begin{equation} \label{qBqJK-rep1}
\widetilde P_n^{A,B}(x;a,b,c;q)=P_n(x;a,b,c;q)-A_{n-1} P_{n-1}
(x;a,b,c;q)-B_{n-1}(x)P_{n-2}(x;aq,bq,cq;q),
\end{equation}
where
$$
A_n=\frac{(abq^2;q)_{n}}{(q;q)_{n}}
\left(A\widetilde{P}_{n+1}^{A,B}(cq)\frac{(cq;q)_{n}}
{c^{n}(abc^{-1}q;q)_{n}}+B\widetilde{P}_{n+1}^{A,B}(aq)
\frac{(aq;q)_{n}}{a^{n}(bq;q)_{n}}\right),
$$
and
$$
B_n(x)=\frac{(abq^2;q)_n(1-q^n)}{(1-aq)(1-cq)(q;q)_nq^n}
\left(A\widetilde{P}_{n+1}^{A,B}(cq)\frac{(cq;q)_n(x-aq)}
{c^n(abc^{-1}q;q)_n}+B\widetilde{P}_n^{A,B}(aq)
\frac{(aq;q)_n(x-cq)}{a^n(bq;q)_n}\right).
$$
Before analyzing the following particular case let us
show that these polynomials can be written as a basic
hypergeometric series.
In fact, by definition of the big $q$-Jacobi polynomials
and \eqref{qBqJK-rep1} we obtain
\begin{equation} \label{sum1} \begin{split}
\displaystyle  \widetilde  P_n^{A,B}(x;a,b,c;q) =& \displaystyle
\sum_{k=0}^{\infty}\frac{(q^{-n},abq^{n},x;q)_kq^k}{(aq^2,cq^2,q;
q)_k}\left(\frac{(1-aq^{k+1})(1-cq^{k+1})}{(1-aq)(1-cq)
}\left[\frac{1-abq^{n+k}}{1-abq^n}\right.\right.\\[0.4cm]&
\displaystyle \left.\left.+A_n\frac{1-q^{-n+k}}{1-q^{-n}}
\right]+B_n(x)\frac{(1-q^{-n+k})(1-q^{-n+k+1})}{(1-q^{-n})(1-
q^{-n+1})}\right).\end{split}
\end{equation}
Now, if we use the identity $(q^{\alpha+1};q)_m (1-q^{\alpha})=
(q^{\alpha};q)_m (1-q^{\alpha+m})$ as well as the fact that the
polynomial on $q^k$ at the RHS, namely $\pi_3$, has three zeros,
namely $q^{\alpha_1}$, $q^{\alpha_2}$, and $q^{\alpha_3}$
which depend, in general, of all parameters, i.e.
$\alpha_{1,2,3}:=\alpha_{1,2,3}(n,x;a,b,c,A,B;q)$, and
$$
\pi_3(q^k)=r(x)(q^k-q^{\alpha_1})(q^k-q^{\alpha_2})
(q^k-q^{\alpha_3}), \qquad \deg r(x)=1.
$$
Then we get
\begin{equation} \label{qBqJK-rep-bhs}
\widetilde  P_n^{A,B}(x;a,b,c;q)=
\widetilde D_n(x){}_6 \varphi_5\left.\left(\begin{array}{c}q^{-n}
, abq^{n},  q^{1-\alpha_1}, q^{1-\alpha_2}, q^{1-\alpha_3},  x \\
aq^2, cq^2, q^{-\alpha_1}, q^{-\alpha_2}, q^{1-\alpha_3}
\end{array}\right|\ q;q\right),
\end{equation}
where
$$
\widetilde D_n(x)=(1-q^{\alpha_1})(1-q^{\alpha_2})(1-q^{\alpha_3})
r(x).
$$
\bn
Notice that $\widetilde D_n$ is, in general, a polynomial of
degree 1 in $x$.
To see that $\widetilde P^{A,B}_n$ is a polynomial of degree $n$
we only need to evaluate \eqref{sum1} for $k=n$ since in this
case the second term on the last bracket vanishes.
\en
%
%  CASO 2
%
\noindent {\bf 2.} The case $x_0=aq$ and $x_1=1$. For these values,
$$
\widehat P^{A,B}_n(x;a,b,c;q):=P_n(x;a,b,c;q)-A\widehat
P^{A,B}_n(aq)K^{BqJ}_{n-1}(x,aq)- B\widehat P^{A,B}_n(1)
K^{BqJ}_{n-1}(x,1).
$$
Then, using  \eqref{BqJ-3} and the relation \eqref{rel-for-BqJ},
we get
{\changesran{
\begin{equation} \label{qBqJ2-rep1}
\begin{array}{rl}
\widehat  P_n^{A,B}(x;a,b,c;q)=& \displaystyle P_{n}(x;a,b,c;q)-
\widehat A_{n-1}P_{n-1}(x;a,b,c;q)\\ &-\widehat B_{n-1}(x)
P_{n-2}(x;aq,bq,cq;q)-\widehat C_{n-1}(x)P_{n-2}(qx;aq,bq,cq;q)
\end{array} \end{equation}}}
where
\begin{align*}
\widehat A_n=&\displaystyle \frac{(aq,abq^2;q)_n}{(q,bq;q)_na^n}
\left(A\widehat P_{n+1}^{A,B}(aq)+\frac{(-1)^n(cq;q)_n}
{(abc^{-1}q;q)_n c^n q^{n^2/2+3n/2}}B\widehat P_{n+1}^{A,B}(1)\right), \\
\widehat B_n(x)=& \displaystyle \frac{(aq,abq^2;q)_n(1-q^n)(x-cq)}
{(q,bq;q)_n(aq)^n(1-aq)(1-cq)}A\widehat P_{n+1}^{A,B}(aq),\\
\widehat C_n(x)=& \displaystyle \frac{(-1)^n(aq,abq^2,cq;q)_n
(1-q^n)aq^3(bx-c)}{(q,bq,abc^{-1}q;q)_na^n c^nq^{n^2/2+3n/2}(1-aq)
(1-cq)}B\widehat P_{n+1}^{A,B}(1).
\end{align*}
Hence, following the same idea of the previous case we get, this
family admits another representation in terms of basic
hypergeometric series.
Indeed some calculations, in the same fashion as in the previous
case, yield
\begin{equation} \label{qBqJ2-rep-bhs}
\widehat P_n^{A,B}(x;a,b,c;q)=\widehat D_n(x){}_6 \varphi_5\left.
\left(\begin{array}{c}q^{-n}, abq^{n}, q^{1-\beta_1},
q^{1-\beta_2}, q^{1-\beta_3}, x \\ aq^2, cq^2, q^{-\beta_1},
q^{-\beta_2}, q^{-\beta_3} \end{array}\right| \ q;q\right),
\end{equation}
where $\widehat D_n(x)$ depends of the parameters defined for this
family, and $q^{\beta_1}$, $q^{\beta_2}$, and $q^{\beta_3}$ are
the zeros of a certain cubic polynomial on $q^k$,
$\beta_i:=\beta_i(n,x;a,b,c,A,B;q)$, $i=1,2,3$, obtained as before
from the expression \eqref{qBqJ2-rep1} and the basic series
representation of the big $q$-Jacobi polynomials.
}

\bn
Observe that, although $\widehat D_n$ is a polynomial of
degree 2 in $x$, $\widetilde P^{A,B}_n$ is a polynomial
of degree $n$ in $x$ since the evaluation of the basic
hypergeometric series \eqref{qBqJ2-rep-bhs} for both
$k=n-1$ and $k=n$ vanishes.
\en

\noindent Two particular interesting cases are the following:
Setting $A=0$ in the $q$-Koornwinder polynomials
\eqref{qBqJK-rep-bhs} we obtain
$$
\widetilde P^B(x;a,b,c;q)=\widetilde D^B_n(x){}_6 \varphi_5 \left.
\left(\begin{array}{c}q^{-n}, abq^{n}, q^{1-\alpha_1},
q^{1-\alpha_2}, q^{1-\alpha_3}, x \\ aq^2, cq^2, q^{-\alpha_1},
q^{-\beta_2}, q^{-\beta_3} \end{array}\right|\ q;q\right),
$$
and setting $A=0$ in the second family \eqref{qBqJ2-rep-bhs} we get
$$
\widehat P^B(x;a,b,c;q)=\widehat D^B_n(x){}_6 \varphi_5 \left.
\left(\begin{array}{c}q^{-n}, abq^{n}, q^{1-\beta_1},
q^{1-\beta_2}, q^{1-\beta_3}, x \\ aq^2, cq^2, q^{-\beta_1},
q^{-\beta_2}, q^{-\beta_3}\end{array}\right| \ q;q\right).
$$

Setting in all the above formulas $c=q^{-N-1}$ we obtain the
$q$-Hahn-Krall polynomials.

Before continuing let us point out that  the above families
satisfy a three-term recurrence relation and a second order linear
difference equation. For more details see \cite{ran-pet}.

%%%%%%%%%%%%%%%%%%%%%%%%%%%%%%%%%%%%%%%%%%%%%%%%%%%%%%%%%%%%%%%%%%%
\subsection{Examples adding one mass point}
\subsubsection{The big $q$-Laguerre-Krall polynomials}
It is a particular case of the $q$-Krall big $q$-Jacobi.
In this case the linear functional ${\cal U}^{BqL}$ is
$$
\pe {{\cal U}^{BqL}} {P}= \pe {{\cal C}^{BqL}} {P}
+ A P(aq), \quad A\ge 0,
$$
where ${\cal C}^{BqL}$ is the functional with respect the big
$q$-Laguerre are orthogonal.
The explicit expression for the polynomials is
$$
\widetilde P_n^A(x;a,c;q)_q=P_n(x;a,c;q)_q-A\frac{P_n(aq;a,c;q)_q
K_{n-1}^{BqL}(x,aq)} {1+A\mathbb{K}^{BqL}_{n-1}(aq)},
$$
or, equivalently, putting $B=0$ in \eqref{qBqJK-rep1} and set
$b=0$,
{\color{blue} 
\[
%\begin{equation}\label{p-BqL}
\begin{split}
\widetilde P_n^{A}(x;a,c;q) = & P_n(x;a,c;q)_q-\frac{A\big(-c
q^{\racion{n+1}2}\big)^n(aq;q)_{n-1}}{1+A\mathbb{K}^{BqL}_{n-1}(aq)}\\
& \times\left[\frac{P_{n-1}(x;a,c;q)_q-(x-cq)(1-q^{-1}) {\cal
D}_{q} P_{n-1}(q^{-1}x;a,c;q)_q}{(cq;q)_{n}(q;q)_{n-1}a^{n-1}}\right]
\end{split}
\]}%\end{equation}
They can be represented as a ${}_6\varphi_5$ basic series.

\subsubsection{The little $q$-Jacobi-Krall polynomials}
These polynomials are orthogonal with respect to the linear
functional ${\cal U}^{lqJ}$
$$
\pe {{\cal U}^{lqJ}} {P}= \pe {{\cal C}^{lqJ}} {P} + AP(0),
%=\sum_{k=0}^{\infty} P(s)\rho(s)q^s \kappa_q+AP(0),
\quad A\ge 0,
$$
where ${\cal C}^{lqJ}$ is the functional of little $q$-Jacobi
polynomials.
The representation formulas for this family is (see \eqref{FR-ker},
\eqref{lqJ})
{\color{blue} \[
\begin{split}
\widetilde p_n^A(x;a,c|q)_q = & p_n(x;a,c|q)_q-A\frac{p_n(0;a,c|q)_q
K_{n-1}^{lqJ}(x,0)} {1+A\mathbb{K}^{lqJ}_{n-1}(0)}= p_n(x;a,c|q)_q
- \frac{A(aq,abq;q)_{n-1}} {1+A\mathbb{K}^{lqJ}_{n-1}(0)}\\ &
\times \left[\frac{(1-abq^{n})P_{n-1}(x;a,b|q)_q-(x-1)(1-q^{-1})
{\cal D}_q P_{n-1}(q^{-1}x;a,b|q)_q}{(1-abq)(q,bq;q)_{n-1} a^{n-1}}\right].
\end{split}
\]
This case leads to a ${}_6 \varphi_5$ basic series.}
\subsubsection{The $q$-Meixner-Krall polynomials}
These polynomials are orthogonal with respect to the linear
functional ${\cal U}^{qM}$
$$
\pe {{\cal U}^{qM}} {P}= \pe {{\cal C}^{qM}} {P} + A P(1),
%=-\sum_{k=0}^{\infty}P(s)\rho(s)q^{-s}\kappa_q+A P(1),
\quad A\ge 0,
$$
where ${\cal C}^{qM}$ is the functional of the $q$-Meixner
polynomials.
The explicit expression for this family is (see \eqref{FR-ker},
\eqref{M})
{\color{blue} \[
\begin{split}
\widetilde M_n^A(x;a,b;q)_q=& M_n(x;a,b|q)_q-A\frac{M_n(1;a,b|q)_q
K_{n-1}^{qM}(x,1)} {1+A\mathbb{K}^{qM}_{n-1}(1)}=M_n(x;a,b|q)_q-
\frac{A(bq;q)_{n-1}}{1+A\mathbb{K}^{qM}_{n-1}(1)} \\ & \times\left[
\frac{M_{n-1}(x;b,c;q)-(x+bc)(1-q){\cal D}_{q^{-1}} M_{n-1}(x;b,c;q)}
{(q,-c^{-1}q;q)_{n-1}}\right]. \end{split}
\]
%\end{small}%
And, this case leads to a ${}_5 \varphi_4$  basic series.}
\subsubsection{The Al-Salam-Carlitz-Krall I polynomials}
These polynomials are orthogonal with respect
to the linear functional ${\cal U}^{ACI}$
$$
\pe {{\cal U}^{ACI}} {P}= \pe {{\cal C}^{ACI}} { P} + A P(1),
\quad A\ge 0,
$$
where ${\cal C}^{ACI}$ is the functional of the Al Salam Carlitz
I polynomials.
The representation formula for this family is (see \eqref{FR-ker},
\eqref{ACI})
{\color{blue} \[
\begin{split}
\widetilde U_n^{(a),A}(x;q)_q=& U_n^{(a)}(x;q)_q -A\frac{U_n^{(a)}
(1;q)_q K_{n-1}^{ACI}(x,1)} {1+A\mathbb{K}^{ACI}_{n-1}(1)}=
U_n^{(a)}(x;q)_q-A\frac{U_n^{(a)}(1;q)_q}{1+A\mathbb{K}^{ACI}_{n-1}
(1)} \\ & \times\frac{q^{n-1}}{(q;q)_{n-1}}\left[U_{n-1}^{(a)}(x;q)-(x-a)
(1-q^{-1}){\cal D}_q U_{n-1}^{(a)}(q^{-1}x;q)\right]. \end{split}
\]
This case leads to a ${}_5\varphi_4$ basic series.}
This family was considered in \cite{ran-pet}.
Since the Al-Salam-Carlitz II are related with the
Al-Salam-Carlitz I by the change $q\to q^{-1}$ the corresponding
$q$-Krall family can be obtained by the same change.

\subsubsection{The little $q$-Laguerre-Krall/Wall-Krall
polynomials}\label{wall-k}
These polynomials are orthogonal with respect to the linear
functional ${\cal U}^{lqL}$
$$
\pe {{\cal U}^{lqL}} {P}= \pe {{\cal C}^{lqL}} { P} + A P(0),
\quad A\ge 0,
$$
where ${\cal C}^{lqL}$ is the functional of the $q$-Laguerre/Wall
polynomials.
The explicit expression for this family is (see \eqref{FR-ker},
\eqref{lqL})
\changeOK{\[
\begin{array}{rl}
\widetilde p_n^{A}(x;a|q)_q = & \displaystyle p_n(x;a|q)_q-\frac
{p_n(0;a|q)_q K_{n-1}^{lqL}(x,0)} {1+A\mathbb{K}^{lqL}_{n-1}(0)}
=\displaystyle p_n(x;a|q)_q-\frac{A}{1+A
\mathbb{K}^{lqL}_{n-1}(0)}\\[4mm]\times & \dfrac{ (aq;q)_{n-1}}{(q;q)_{n-1}
a^{n-1}q^{n-1}}\left[p_{n-1}(x;a|q)
-a(1-q)q^{n-1}{\cal D}_q p_{n-1}(x;a|q)\right].
\end{array}
\]
This case leads to a ${}_4\varphi_3$ basic series.}
\subsubsection{The $q$-Laguerre-Krall polynomials}
These polynomials are orthogonal with respect to the linear
functional ${\cal U}^{qL}$
$$
\pe {{\cal U}^{qL}} {P}= \pe {{\cal C}^{qL}} {P} + AP(0),
\quad A\ge 0,
$$
where ${\cal C}^{qL}$ is the functional of the $q$-Laguerre
polynomials. In this case \eqref{FR-ker} and  \eqref{L} yield
{\color{blue} \[
\begin{array}{rl}
\widetilde L_n^{(\alpha),A}(x;q)_q=& \displaystyle
L_n^{(\alpha)}(x;q)_q-A\frac{L_n^{(\alpha)}(0;q)_q K_{n-1}^{qL}
(x,0)} {1+A\mathbb{K}^{qL}_{n-1}(0)} \\[4mm] = & \displaystyle
L_n^{(\alpha)}(x;q)_q - A\frac{L_n^{(\alpha)}(0;q)_q}
{1+A\mathbb{K}^{qL}_{n-1}(0)}\left[q^{n-1}L_{n-1}^{(\alpha)}(x;q)
-\frac{q^{-1}-1} {a}{\cal D}_q L_{n-1}^{(\alpha)}(q^{-1}x;q)\right].
\end{array}
\]
This case leads to a ${}_4\varphi_3$ basic series.}
\subsubsection{The $q$-Charlier-Krall polynomials}
These polynomials are orthogonal with respect to the linear
functional ${\cal U}^{qC}$
$$
\pe {{\cal U}^{qC}} {P}= \pe {{\cal C}^{qC}} {P} + A P(1),
\quad A\ge 0,
$$
where ${\cal C}^{qC}$ is the weight function of the
$q$-Charlier polynomials.  For these polynomials \eqref{FR-ker}
and  \eqref{C} yield
{\color{blue} \[
\begin{array}{rl}
\widetilde C_n^{A}(x;a;q)_q=& \displaystyle C_n(x;a;q)_q-A
\frac{C_n(1;a;q)_qK_{n-1}^{qC}(x,1)} {1+A\mathbb{K}^{qC}_{n-1}(1)}
\\[4mm] = & \displaystyle C_n(x;a;q)_q-\frac{A}
{1+A\mathbb{K}^{qC}_{n-1}(1)}\left[\frac{C_{n-1}(x;a;q)-x(1-q)
{\cal D}_{q^{-1}} C_{n-1}(x;a;q)}{(-a^{-1}q,q;q)_{n-1}}\right].
\end{array}
\]
This case leads to a ${}_4 \varphi_3$ basic series.}
\subsubsection{The Stieltjes-Wigert-Krall polynomials}\label{stiWie}
These polynomials are orthogonal with respect to the linear
functional ${\cal U}^{SW}$
$$
\pe {{\cal U}^{SW}} {P}= \pe {{\cal C}^{SW}} { P} + A P(0),
\quad A\ge 0,
$$
where ${\cal C}^{SW}$ is the functional of the Stieltjes-Wigert
polynomials. The representation formula for this family has the
form  (see \eqref{FR-ker}, \eqref{SW})
{\color{blue} 
\[
\widetilde S_n^{A}(x;q)_q = %S_n(x;q)_q - A\frac{S_n(0;q)_q
%K_{n-1}^{SW}(x,0)} {1+A\mathbb{K}^{SW}_{n-1}(0)}=
S_n(x;q)_q-A\frac{S_n(0;q)_q} {1+A\mathbb{K}^{SW}_{n-1}(0)}\left[
q^{n-1}S_{n-1}(x;q)-(q^{-1}-1){\cal D}_qS_{n-1}(q^{-1}x;q)\right].
\]
This case leads to a ${}_3 \varphi_3$ basic series. }
This family was firstly studied in \cite{chi1}.
%%%%%%%%%%%%%%%%%%%%%%%%%%%%%%%%%%%%%%%%%%%%%%%%%%%%%%%%%%%%%%%%%%%
\subsection{Some algebraic properties of $\widetilde{P}_n^A(s)_q$}
{\color{blue} 
In \cite{ran-pet} it is shown that  the $q$-Krall-type orthogonal
polynomials satisfy a second order linear difference equation of
the form

\begin{equation} \label{q-diff-equiv}
\widetilde{\sigma}(s;n) \widetilde{P}_n(s-1)_q -
\widetilde{\varphi}(s;n)\widetilde{P}_n (s)_q +
\widetilde{\varsigma}(s;n)\widetilde{P}_n (s+1)_q =0,
\end{equation}
where
\begin{equation} \label{q-diff-coef}
\begin{array}{l}
\widetilde{\sigma}(s;n) = t(s;n)  [a(s;n)d(s;n)-c(s;n)b(s;n)],\\[2mm]
\widetilde{\varphi}(s;n)= -\pi(s;n) [c(s;n)f(s;n) - e(s;n)d(s;n)], \\[2mm]
\widetilde{\varsigma}(s;n)= r(s;n) [e(s;n)b(s;n)-a(s;n)f(s;n)],
\end{array}
\end{equation}
being
$$
\begin{array}{l}
r(s;n)=\varsigma(s+1;n)\pi(s+1;n),\quad
c(s;n)= -\sigma(s+1;n)\, b(s+1;n) ,\\[3mm]
d(s;n)=a(s+1;n) \varsigma(s+1;n) + b(s+1;n) \varphi(s+1;n),
\quad t(s;n)=\sigma(s;n)\pi(s-1;n),\\[3mm]
e(s;n)= \sigma(s;n)b(s-1;n) +a(s-1;n)\varphi(s;n),\quad
f(s;n)= -a(s-1;n) \varsigma(s;n),
\end{array}
$$
where $a(s;n)$, $b(s;n)$, and $\pi(s;n)$ are the coefficients of the
representation formula for the Krall-type polynomials  $\widetilde{P}_n$
\begin{equation*}
\pi(s;n) \widetilde{P}_n(s)_q = a(s;n) P_n(s)_q + b(s;n) P_n(s+1)_q,
%\label{rep-for}
\end{equation*}
and $\sigma$, $\varphi$, and $\varsigma$ are the coefficients of the
second order difference equation that the starting polynomials $P_n$ satisfy
\begin{equation*}
\sigma(s;n)P_n(s-1)_q-\varphi(s;n) P_n(s)_q+\varsigma(s;n)P_n(s+1)_q=0.
%\label{eqdif-dis}
\end{equation*}

Also in \cite{ran-pet}  the TTRR for the polynomials
$\widetilde{P}_n^A(s)_q$ is computed
\begin{equation}\label{RRTT}
x(s)\widetilde{P}_n^A(s)_q = \alpha^A_n\widetilde{P}_{n+1}^A
(s)_q + \beta^A_n \widetilde{P}_n^A(s)_q+ \gamma^A_n
\widetilde{P}_{n-1}^A(s)_q,
\end{equation}
where the coefficients $\alpha^A_n, \, \beta^A_n$, and $\gamma^A_n$
are given by
$$
\begin{array}{l}
\displaystyle \alpha^A_n=\alpha_n, \qquad  \gamma^A_n=\alpha_{n-1}
\frac{{\tilde d}_{n}^2}{{\tilde d}_{n-1}^2}\ne 0, \\[6mm]
\displaystyle \beta^A_n = \beta_n +\frac{A P_n(s_0)_q}{d_n^2}
\left(\alpha_n\frac{P_{n+1}(s_0)_q} {1+A\mathbb{K}_n(s_0)}-
\gamma_n\frac{P_{n-1}(s_0)_q}{1+A\mathbb{K}_{n-1}(s_0)}\right),
\end{array}
$$
where $\alpha_n$, $\beta_n$, and $\gamma_n$ are the coefficients
of the TTRR of the starting family of $q$-polynomials
\eqref{q-TTRR},  $d_n^2=\pe{{\cal C}}{P_nP_n}$ and
$$
{\tilde d}^2_n = \langle {\cal U}, \widetilde{P}_n^A
\widetilde{P}_n^A \rangle = d_n^2+\left[A\widetilde{P}^A_n(s_0)_q
\right]^2 \mathbb{K}_{n-1}(s_0)+A\left[\widetilde{P}^A_n(s_0)_q
\right]^2 =\frac{1+A\mathbb{K}_n(s_0)}{1+A\mathbb{K}_{n-1}(s_0)}
d_n^2,
$$
are the square of the norms of the polynomials $P_n$ and
$\widetilde P_n^A$, respectively.

%%%%%%%%%%%%%%%%%%%%%%%%%%%%%%%%%%%%%%%%%%%%%%%%%%%%%%%%%%%%%%%%%%%
\subsubsection{Some examples}
Here we will restrict ourselves to the more simple cases.
The other cases are analogously and we will omit them here.

\bigskip
\noindent\textbf{Little $q$-Laguerre-Krall / Wall-Krall polynomials:}
The Wall polynomials satisfy the following SODE (recall that $x:=q^ s$)
$$
a p_n(qx;a|q)+(1-x)p_n(x/q;a|q)+[(x-a-1)+q^{-n}(1-q^n)x]p_n(x;a|q)=0
$$
as well as the relation \cite[Eq. (3.2), page 175]{ran02}
$$
x{\cal D}_{q} p_n(x;a|q)= \frac{1-q^n}{1-q}\left(p_n(x;a|q)-p_{n-1}(x;a|q)\right)
\Rightarrow (1-q^n)p_{n-1}(x;a|q)=-q^n p_n(x;a|q)+ p_n(qx;a|q).
$$
Then, for the kernel $K_{n-1}(x,0)$ we have the expression
$$
 x K_{n-1}(x,0)=-\frac{(aq;q)_{n+1}}{(1-q^n)(q;q)_n a^n}\left[p_{n}(x;a|q) - p_{n}(qx;a|q)\right].
$$
Thus from the expression for the Wall-Krall polynomials in section \ref{wall-k} we obtain
$$
x \widetilde p_n^{A}(x;a|q)_q = (x+ b_n^A) \displaystyle p_n(x;a|q)_q - b_n^A p_n(qx;a|q)_q,
$$
\changeOK{
$$x\widetilde p_n^{A}(x;a|q)_q=(x-\frac{c_n}{q^{n-1}})p_n(x;a|q)+
ac_n(\frac{1-q^{n-1}}{q^{n-1}}-\frac{1+a(1-q^{n-1}-q^n)}{a})
\frac{(q^np_n(x;a|q)-p_n(qx;a|q))}{(1-q^n)(1-aq^n)}$$
where 
$$
c_n=\frac{A(aq;q)_{n-1}}{(1+A\mathbb K_{n-1}(0))a^{n-1}(q;q)_{n-1}}.
$$}

where
$$
b_n^A= A\frac {p_n(0;a|q)_q} {1+A\mathbb{K}^{lqL}_{n-1}(0)} \frac{(aq;q)_{n+1}}{(1-q^n)(q;q)_n a^n}.
$$
Therefore, Theorem 2 of \cite[page 60-61]{ran-pet} gives the
following values for the coefficients of the SODE \eqref{q-diff-equiv}
\begin{small}\[
\begin{split}
\widetilde{\varsigma}(s;n)= & -\frac{a\left(x+1\right)\left(q^{n}x (b_n^A)^2-x(b_n^A)^2-x^2 b_n^A-a q^{n} xb_n^A+q^{n}xb_n^A+xb_n^A-q^{n}b_n^A-aq^{n}x^2+aq^{n}x\right)}
{q^{n}},\\
\widetilde{\varphi}(s;n)= & \,
x \left[\left(\left(1\!-\!q^{-n}\right)x\!-\!x+a+1 \right)\left(b_n^A+x\!-\!1\right)\!-\!\left(1\!-\!x\right)
b_n^A\right]
\left(a \left(b_n^A+x+1\right)\!-\!\left(\left(1\!-\!q^{-n}\right)\left(x+ 1\right)\!-\!x+a\right)b_n^A\right)\\
& +ax^2 \left(-b_n^A-x+1\right)b_n^A,\\
\widetilde{\sigma}(s;n)= &
\frac{\left(x\!-\!1\right)^2\left(q^{n}x(b_n^A)^2\!-\!x (b_n^A)^2\!+q^{n}(b_n^A)^2\!-(b_n^A)^2\!-x^2
b_n^A\!-\!aq^{n}xb_n^A+q^{n}xb_n^A\!-\!xb_n^A\!-\!aq^{n}b_n^A\!-\!aq^{n}x^2\!-\!aq^{ n}x\right)}
{q^{n}}.
\end{split}
\]
\end{small}%
For this family we have the following coefficients of the TTRR \eqref{RRTT}
\[
\begin{split}
\alpha^A_n= & -q^n(1-aq^{n+1}), \\
\beta^A_n= & q^n(1-aq^{n+1})+aq^n(1-q^n)) \\ & -A\frac{q^n(1-a
q^{n+1})(1-aq)}{(1-aq)d_n^2+A(1-aq^{n+1})}+A\frac{q^{n-1}(1-a
q^{n})(1-aq)}{(1-aq)d_{n-1}^2+A(1-aq^{n})},\\
\gamma^A_n= &
\alpha_{n-1}^A \frac{(d_n^2\big((1-aq)d_{n-1}^2+1-aq^n\big)+
(1-aq)d_{n-1}^2)((1-aq)d_{n-2}^2+1-aq^{n-1})}{(d_{n-1}^2
\big((1-aq)d_{n-2}^2+1-aq^{n-1}\big)+(1-aq)d_{n-2}^2)((1-aq)
d_{n-1}^2+1-aq^{n})},
\end{split}
\]
where $d_n^2=\frac{(q;q)_n}{(aq;q)_n} (aq)^n$.

\bigskip
\noindent\textbf{$q$-Stieltjes-Wigert-Krall polynomials:}
The $q$-Stieltjes-Wigert  polynomials satisfy the SODE
$$
x S_n(qx;q)+  S_n(x/q;q) -[1+q^n x] S_n(x;q)=0.
$$
Now, combining relation \cite[Eq. (3.2), page 175]{ran02} and
the TTRR of the Stieltjes-Wigert polynomials we have
$$
x^2{\cal D}_{q} S_n(x;q)= \frac{1-q^n}{1-q}(q^n x+q^2)q^{-n}s_n(x;q)-
\frac{q^{-n}}{1-q}S_{n-1}(x;q)
\Rightarrow
$$
$$
S_{n-1}(x;q)=  - [q^{2n}x+q^2(q^n-1)] S_n(x;q)-x q^n S_n(q x;q).
$$
Then, for the kernel $K_{n-1}(x,0)$ we have the expression
$$
 x K_{n-1}(x,0)=-(1-q^{n+1})q^{-3n-1} \left[ \big((1-q^n)(1-q^2)+q^{2n}x\big) S_n(x;q)-
x q^n S_n(q x;q)  \right].
$$
Thus from the expression for the $q$-Stieltjes-Wigert-Krall polynomials in section \ref{stiWie} we obtain
$$
x \widetilde S_n^{A}(x;q)_q =\left[ (1+ b_n^Aq^{2n})x +(1-q^n)(1-q^2) b_n^A \right]
S_n(x;a)_q + b_n^A q^n x S_n(qx;q)_q,
$$
being
$$
b_n^A=  Aq^{-3n-1}\frac{S_n(0;q)_q(1-q^{n+1})}{1+A\mathbb{K}^{SW}_{n-1}(0)}.
$$
Therefore, Theorem 2 of \cite[page 60-61]{ran-pet} gives the
following values \eqref{q-diff-coef} for the coefficients of the SODE \eqref{q-diff-equiv}
which explicit expression we will omit here.

For this family we have the following coefficients of the TTRR \eqref{RRTT}

\[
\begin{split}
\alpha^A_n= & -(1-q^{n+1})q^{-2n-1}, \\
\beta^A_n= & (1+q-q^{n+1})q^{-2n-1} -\frac{A}{q^{3n+1}}\left(
\frac{1}{(q;q)_{n}+A}-\frac{q^3}{(q;q)_{n-1}+A}\right),\\
\gamma^A_n= & -\frac{1}{q^{2n+2}}\frac{((q;q)_n+A)((q;q)_{n-1}
+A(1-q^{n-1}))}{((q;q)_{n-1}+A)((q;q)_n+A(1-q^{n}))}.
\end{split}
\]
}
%%%%%%%%%%%%%%%%%%%%%%%%%%%%%%%%%%%%%%%%%
\section{Limit relations between $q$-Krall-type orthogonal
polynomials}
\noindent
In this section, we study the limit relations involving the
$q$-Krall-type orthogonal polynomials associated with some families
of $q$-polynomials of the $q$-Hahn Tableau \cite{koo94,ran02}.
As we already pointed out the
$q$-Koornwinder polynomials $\widetilde P^{A,B}_n(x;a,b,c;q)$
\eqref{qBqJK-rep1} is the $q$-analogue of the Koornwinder
polynomials $P^{A,B}_n(x)$
\cite{koo84}. In fact, a direct calculation show
$$
\lim_{q\to1-}\widetilde P^{A,B}_n(x;a,b,c;q)= P^{A,B}_n(x).
$$
Let now consider the other limits.
\be
\item {\bf  Big $q$-Jacobi $\to$ Big $q$-Laguerre.}
We know that the big $q$-Laguerre is a special case of big
$q$-Jacobi setting $b=0$, i.e. $P_n(x;a,0,c;q)=P_n(x;a,c;q)$.
Then, from \eqref{pol_one_mas} we get
$$
\widetilde P^A_n(x;a,0,c;q)=\widetilde P^A_n(x;a,c;q).
$$
\item {\bf  Big $q$-Jacobi $\to$ Little $q$-Jacobi}.
The little $q$-Jacobi polynomials can be obtained from the big
$q$-Jacobi polynomials by linear change of the variable
$x\to cqx$ and taking the limit $c\to \infty$, i.e.
$\lim_{c\to \infty} P_n(cqx;a,b,c;q)=p_n(x;a,b|q)$.
In this case, putting $xcq=aq$ and taking the limit $c\to \infty$
we get $x\to 0$, thus $\displaystyle \lim_{c\to \infty}
P_n(aq;a,b,c;q)$ $=p_n(0;a,b|q)$. Taking into account that the
the norm of big $q$-Jacobi transforms into the norm of the little
$q$-Jacobi we obtain
$$
\lim_{c\to \infty} \widetilde P^A_n(cqx;a,b,c;q)=\widetilde
p^A_n(x;a,b|q).
$$
\item {\bf  Big $q$-Jacobi $\to$ $q$-Meixner.}
If we take the limit $a\to \infty$ in the big $q$-Jacobi we obtain
the $q$-Meixner polynomials \cite{ks}. Thus, from
\eqref{pol_one_mas} we deduce
$$
\lim_{a\to \infty} \widetilde P^A_n(q^{-s};a,b,c;q)=\widetilde
M^A_n(q^{-s};c,-b^{-1};q).
$$
\item {\bf Big $q$-Jacobi $\to$ Hahn}.
Setting $c=q^{-N-1}$ in the big $q$-Jacobi we get the
$q$-Hahn po\-ly\-no\-mials $\widehat P_n^{A,B}(x;a,b,q^{-N-1};q)=
\widehat Q_n^{A,B}(x;a,b,N|q)$.
Substituting $x=q^{-x}$, $a=q^{\alpha}$, $b=q^\beta$, we recover
the Hahn-Krall polynomials studied in \cite{ran-mar2}
$\lim_{q\to 1^-} \widehat Q_n^{0,A}(q^{-x};q^\alpha,q^\beta,
q^{-N-1}|q)$ $=Q_n^{0,A}(x;\alpha,\beta,N)$.
Notice that from the Hahn-Krall polynomials it is possible to
obtain several other families of Krall-type polynomials via
appropriate limits (see \cite{ran-mar3}).
\item {\bf  Big $q$-Laguerre $\to$ Al-Salam-Carlitz I}.
Substituting $x\to aqx$ and $c\to ac$ in the big $q$-Laguerre
polynomials and taking $a\to 0$ we obtain the Al-Salam-Carlitz I
polynomials $\lim_{a\to 0} \frac{P_n(aqx;a,ac;q)}{a^n}=q^n
U_n^{(c)}(x;q)$.
Therefore,
$$
\lim_{a\to 0} \frac{\widetilde P^A_n(aqx;a,ac;q)}{a^n}=q^n
\widetilde U_n^{(c),A}(x;q).
$$
\item {\bf  Big $q$-Laguerre $\to$ Little $q$-Laguerre/Wall}.
The little $q$-Laguerre polynomials can be obtained from the big
$q$-Laguerre polynomials by setting $x\to bqx$ and then taking the
limit $b\to\infty$: $\lim_{b\to \infty} P_n(bqx;a,b;q)=p_n(x;a|q)$.
Thus
$$
\lim_{b\to \infty} \widetilde P^A_n(bqx;a,b;q)=\widetilde
p^A_n(x;a|q).
$$
\item {\bf  Little $q$-Jacobi $\to$ Little $q$-Laguerre/Wall}.
Setting $b=0$ in the little $q$-Jacobi polynomials we get the
little $q$-Laguerre $p_n(x;a,0|q)=p_n(x;a|q)$, then
$$
\widetilde p^A_n(x;a,0|q)= \widetilde p^A_n(x;a|q).
$$
\item {\bf  Little $q$-Jacobi $\to$ $q$-Laguerre}.
In this case straightforward calculations give us
$$
\lim_{b\to \infty} \widetilde p^A_n\left(-\frac x
{bq};q^\alpha,b\Big|\, q\right)= \frac{(q;q)_n}
{(q^{\alpha+1};q)_n}\widetilde L_n^{(\alpha),A}(x;q).
$$
\item {\bf  $q$-Meixner $\to$ $q$-Laguerre}.
Straightforward calculations yield
$$
\lim_{c\to \infty} \widetilde M^A_n(cax;a,c;q)=\frac{(q;q)_n}
{(q^{\alpha+1};q)_n}\widetilde L_n^{(\alpha),A}(x;q).
$$
\item {\bf  $q$-Meixner $\to$ $q$-Charlier}. $
\lim_{b\to 0} \widetilde M^A_n(x;b,a;q)= \widetilde C^A_n(x;a;q).$
\item {\bf  $q$-Laguerre $\to$ Stieltjes-Wigert}.
$\lim_{\alpha \to \infty} \widetilde
L_n^{(\alpha),A}(xq^{-\alpha};q)=\widetilde S^A_n(x;q).$
\item {\bf  $q$-Charlier $\to$ Stieltjes-Wigert}. $
\lim_{a\to \infty} \widetilde C^A_n(ax;a;q)=(q;q)_n\widetilde
S^A_n(x;q).$
\ee
To finish this work let us point out that for the other families
of the $q$-Hahn tableau, i.e., for the $q$-Kravchuk, alternative
$q$-Charlier the same results can be obtained in an analogous way.%
%%%%%%%%%%%%%%%%%%%%%%%%%%%%%%%%%%%%%%%%%%%%%%%%%%%%%%%%%%%%%%%%%%%
 \begin{center}
 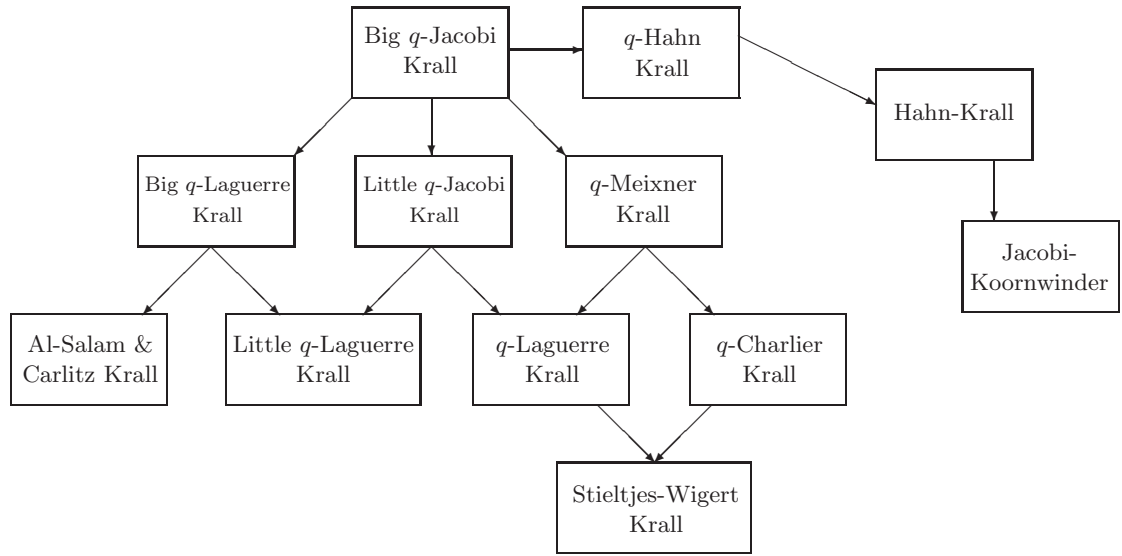
\begin{figure}[t] \setlength{\unitlength}{1900sp}%
\begingroup\makeatletter\ifx\SetFigFont\undefined%
\gdef\SetFigFont#1#2#3#4#5{%
  \reset@font\fontsize{#1}{#2pt}%
  \fontfamily{#3}\fontseries{#4}\fontshape{#5}%
  \selectfont}%
\fi\endgroup%
\begin{picture}(14469,7179)(-731,-7138)
\thinlines
{\put(-719,-5191){\framebox(2025,1170){}} %% ACI
}%
{\put(2091,-5191){\framebox(2525,1170){}}  %%lqL
}%
{\put(5311,-5191){\framebox(2025,1170){}}  %%qL
}%
{\put(3736,-1186){\framebox(2025,1170){}}  %% bqJ
}%
{\put(3736,-1186){\vector(-1,-1){765}} % OK
}%
{\put(946,-3121){\framebox(2025,1170){}}  % BqL
}%
{\put(3781,-3121){\framebox(2025,1170){}} % lqJ
}%
{\put(5761,-1186){\vector(1,-1){765}} % OK
}%
{\put(6526,-3121){\framebox(2025,1170){}} %qM
}%
{\put(4771,-1186){\vector(0,-1){765}}  % OK
}%
{\put(5761,-556){\vector(1, 0){990}} % OK
}%
{\put(6751,-1186){\framebox(2025,1170){}} %qH
}%
{\put(8776,-376){\vector(2,-1){1796}}  % OK
}%
{\put(1891,-3121){\vector(-1,-1){895}}  % OK
}%
{\put(1891,-3121){\vector(1,-1){895}}  % OK
}%
{\put(4771,-3121){\vector(-1,-1){895}} % OK
}%
{\put(4771,-3121){\vector(1,-1){895}} % OK
}%
{\put(7561,-3121){\vector(-1,-1){895}} % OK
}%
{\put(7561,-3121){\vector(1,-1){895}} % OK
}%
{\put(8146,-5191){\framebox(2025,1170){}} % qC
}%
{\put(6406,-7126){\framebox(2525,1170){}} % SW
}%
{\put(6931,-5191){\vector( 1,-1){760}}  % OK
}%
{\put(8428,-5191){\vector(-1,-1){760}} % OK
}%
{\put(10576,-1996){\framebox(2025,1170){}} % OK
}%
{\put(11701,-3976){\framebox(2025,1170){}} % OK
}%
{\put(12106,-1996){\vector( 0,-1){810}} % OK
}%
\put(-884,-4500){\makebox(0,0)[lb]
{\smash{{\SetFigFont{9}{14.4}{\rmdefault}{\mddefault}{\updefault}
{\obec{Al-Salam \& \\[-2mm] Carlitz Krall}}% OK
}}}}
\put(10800,-1476){\makebox(0,0)[lb]{\smash{{\SetFigFont{9}{14.4}
{\rmdefault}{\mddefault}{\updefault}{Hahn-Krall}% OK
}}}}
\put(11481,-3301){\makebox(0,0)[lb]{\smash{{\SetFigFont{9}{14.4}
{\rmdefault}{\mddefault}{\updefault}
{\obec{Jacobi-\\[-2mm]Koornwinder}}% OK
}}}}
\put(3576,-466){\makebox(0,0)[lb]{\smash{{\SetFigFont{9}{14.4}
{\rmdefault}{\mddefault}{\updefault}{ \obec{Big $q$-Jacobi\\[-2mm]
Krall}} % OK
}}}}
\put(7066,-500){\makebox(0,0)[lb]{\smash{{\SetFigFont{9}{14.4}
{\rmdefault}{\mddefault}{\updefault}{\obec{$q$-Hahn\\[-2mm]
Krall}}% OK
}}}}
\put(820,-2401){\makebox(0,0)[lb]{\smash{{\SetFigFont{8}{14.4}
{\rmdefault}{\mddefault}{\updefault}{\obec{Big $q$-Laguerre\\[-2mm]
Krall}} % OK
}}}}
\put(3570,-2401){\makebox(0,0)[lb]
{\smash{{\SetFigFont{8}{14.4}{\rmdefault}{\mddefault}{\updefault}
{\obec{Little $q$-Jacobi\\[-2mm] Krall}}% OK
}}}}
\put(6626,-2401){\makebox(0,0)[lb]{\smash{{\SetFigFont{9}{14.4}
{\rmdefault}{\mddefault}{\updefault}{\obec{$q$-Meixner\\[-2mm]
Krall}}%
}}}}
\put(1970,-4500){\makebox(0,0)[lb]{\smash{{\SetFigFont{9}{14.4}
{\rmdefault}{\mddefault}{\updefault}{\obec{Little $q$-Laguerre
\\[-2mm] Krall}}%  OK
}}}}
\put(5396,-4500){\makebox(0,0)[lb]{\smash{{\SetFigFont{9}{14.4}
{\rmdefault}{\mddefault}{\updefault}{\obec{$q$-Laguerre\\[-2mm]
Krall}}% OK
}}}}
\put(8276,-4500){\makebox(0,0)[lb]{\smash{{\SetFigFont{9}{14.4}
{\rmdefault}{\mddefault}{\updefault}{\obec{$q$-Charlier\\[-2mm]
Krall}}% OK
}}}}
\put(6400,-6461){\makebox(0,0)[lb]{\smash{{\SetFigFont{9}{14.4}
{\rmdefault}{\mddefault}{\updefault}{\obec{Stieltjes-Wigert\\[-2mm]
Krall}}%
}}}}
\end{picture}%
  \caption{The $q$-Hahn-Krall Tableau}
 \end{figure}
 \end{center}%
%%%%%%%%%%%%%%%%%%%%%%%%%%%%%%%%%%%%%%%%%%%%%%%%%%%%%%%%%%%%%%%%%%%
\noindent {\bf Acknowledgements:} Discussions with Prof. F.
Marcellan were very fruitful in order to improve this manuscript.
We thank the unknown referee for his remarks and suggestions that
allow us to improve the paper, and also Profs. M.E.H. Ismail and
A. Zhedanov for pointing out the references \cite{uva69} and
\cite{vin01}, respectively.
This work has been partially supported by Direcci\'{o}n General de
Investigaci\'{o}n del Ministerio de Educaci\'{o}n y Ciencia of
Spain under grant BFM2003-06335-C03-01 (RAN), BFM2003-06335-C03-02
(RCS), and the PAI grant  FQM-0262 (RAN).
%%%%%%%%%%%%%%%%%%%%%%%%%%%%%%%%%%%%%%%%%%%%%%%%%%%%%%%%%%%%%%%%%%%

\end{document}